\documentclass[twocolumn]{svjour3}  
\smartqed  

\usepackage[dvipdfmx]{color}
\usepackage{graphicx}
\usepackage[T1]{fontenc}
\usepackage{mathptmx}
\usepackage[scaled]{helvet}
\usepackage{courier}

\makeatletter
\let\cl@chapter\undefined
\makeatletter
\usepackage{amsmath,amssymb}   
\usepackage{mathtools} 
\usepackage{cleveref}  
\Crefname{equation}{Eq.}{Eqs.}%
\Crefname{figure}{Fig.}{Figs.}%
\usepackage{amsbsy}
\usepackage{bm}
\usepackage{algorithmic}
\usepackage{algorithm}

\journalname{Journal of Marine Science and Technology}
\begin{document}

\title{Simultaneous optimization of control gains and reference filter coefficients for trajectory tracking control
%
}
\author{Amane Sakanashi \and Rin Suyama \and
        Atsuo Maki      \and Youhei Akimoto
}

\institute{ A.Sakanashi \and R.Suyama \and  A.Maki \at
              Osaka University, 2-1 Yamadaoka, Suita, Osaka, Japan \\
              \email{sakanashi\_amane@naoe.eng.osaka-u.ac.jp }\\
              \email{maki@naoe.eng.osaka-u.ac.jp}
           \and
           Y.Akimoto \at
           Faculty of Engineering, Information and Systems, University of Tsukuba, 1-1-1 Tennodai, Tsukuba, Ibaraki 305-8573,\\
           RIKEN Center for Advanced Intelligence Project, 1-4-1 Nihonbashi, Chuo-ku, Tokyo 103-0027, Japan\\
           \email{akimoto@cs.tsukuba.ac.jp} 
}

\date{Received: date / Accepted: date}

\maketitle

\begin{abstract}
Research on vessel automation and autonomy is currently being conducted by various countries and institutions. Safe and accurate ship control algorithms are crucial to realize automated operation. Actuator drive constraints of a target ship may jeopardize the stability of the control law and require complex theory. In this study, we include a penalty term to the control law gain optimization stage of dynamic positioning systems to account for the amounts by which the actuator input value and its rate of change exceed the constraint. The parameters for generating a suitable reference path for the control law are identified simultaneously with the control gains. The simulation results show that the proposed method can realize control parameters and a reference design with excellent tracking performance while determining the cost of the controller design by considering the effects of both the actuators and rate saturation.

\end{abstract}

\section{Introduction}
Maneuvering a ship is a highly complex task that requires simultaneous control of multiple actuators, including the propeller, rudder, and thrusters. Automation of ship control has been studied for many years to reduce the burden on crew members and prevent human error. The most commonly used automatic control algorithm in the world is the proportional-integral-derivative (PID) control. Ship autopilot technology was initially built based on the theoretical analysis of PID control. However, the simple mechanism of PID control does not always give satisfactory results when more complicated or high-performance control is oriented \cite{Minorsky1922}. For this reason, various control theories have been developed and employed in this technology, including, for example, those based on the Lyapunov stability theory.
 \par
However, even if a control law has been obtained theoretically, the control gains and other parameters contained within the control law must be tuned before it can be used in practical applications. Tuning control gains by trial and error is the most intuitive method. Such parameter-tuning methods for control laws have been studied since the 1990s. For example, Tomera optimized the parameters of the PID control using an ant colony algorithm \cite{Tomera2014}. Witkowska optimized the parameters of the backstepping control using a genetic algorithm \cite{Witkowska2012}.
 \par
 %
 This study addresses parameter tuning of dynamic positioning system (DPS) control laws for ships. DPS is used to control a vessel to maintain a fixed point and track its path using active thrusters such as propellers. The first DPS was developed in the 1960s. It was a PID control combined with a low-pass filter for the three horizontal degrees of freedom, which could cause phase changes that could affect the system's stability. In the 1970s, Belchen proposed a more advanced control algorithm based on the linear Kalman filter and multivariate linear optimization control theories \cite{Balchen1980}. This method has been developed in various ways \cite{Fung1983}\cite{Sorensen1996}. However, it has drawbacks, such as the dynamics must be linearized under certain conditions, and the stability of the DPS is not guaranteed \cite{Saelid1983}. In the 1990s, dynamic positioning (DP)  nonlinear control was studied extensively following the development of nonlinear control theory. Nonlinearity in DP can reduce the control complexity by removing the complex linearization and its theory. Stephans et al. proposed the Fuzzy control law \cite{Stephens1995} . Araset et al. proposed a nonlinear feedback linearization and backstepping of DP \cite{Aarset1998}. Agostinho et al. proposed the introduction of nonlinear sliding mode control  to DP \cite{Agostinho2009}. Fossen developed DP control by nonlinear observer by ignoring disturbances  \cite{Fossen1998}. The same author proposed a passive nonlinear observer with filtering to estimate the low-frequency position and speed of a vessel under disturbance \cite{Fossen1999}. Loria combined this observer with the PID control law to evaluate the output feedback \cite{Loria2000}.
\par

\begin{figure*}[tb]
        \centering 
        \includegraphics[width=0.7\hsize]{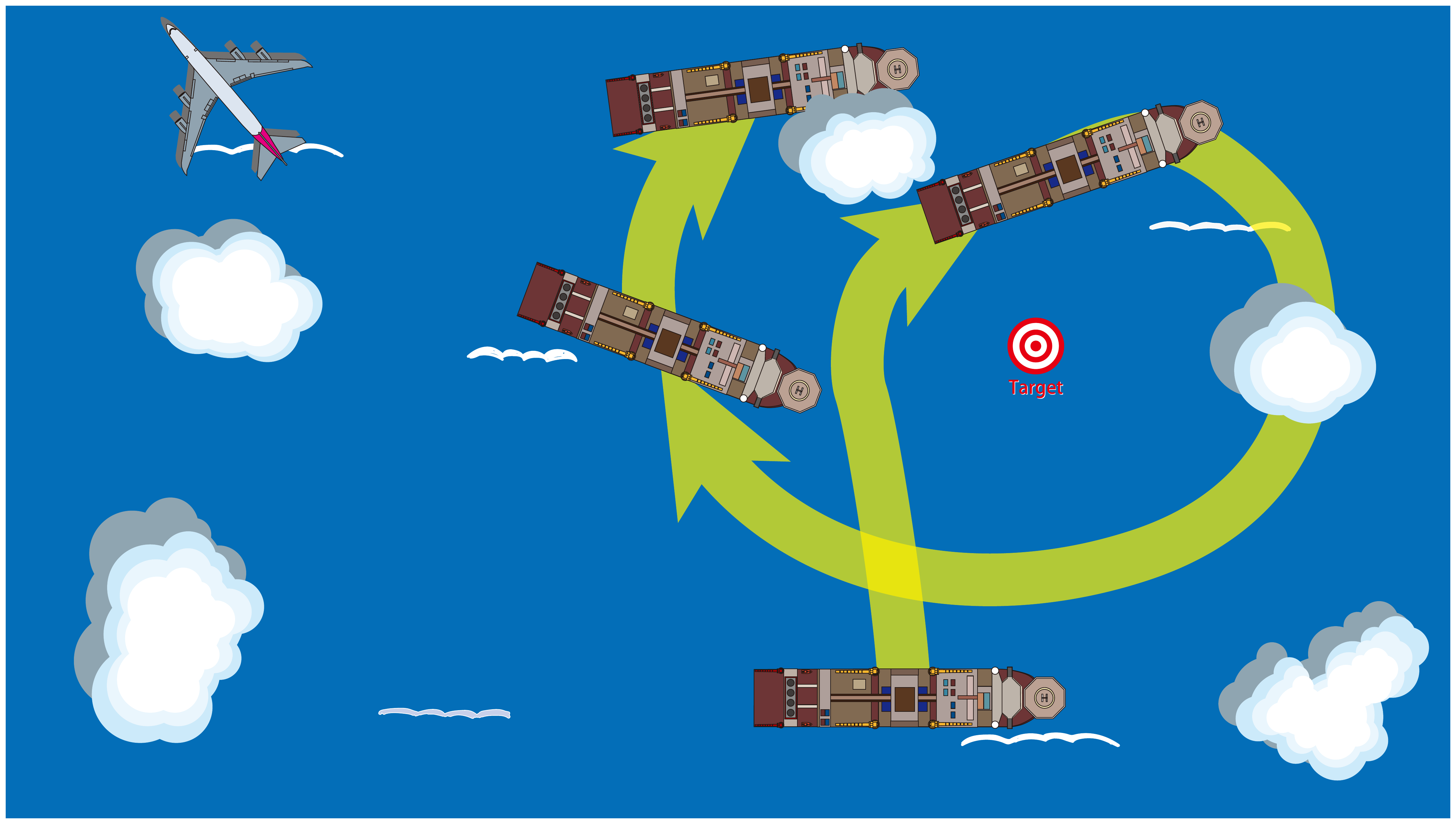}
        \caption{
         Example of unstable controlling of a ship due to actuator saturation
        }
        \label{fig:Example of Actuator Satuation}
    \end{figure*}

%
In the design of such control laws , the treatment of constraints on actuators, such as rudders and propellers, is an issue. In this paper, we refer to the situation where the upper and lower limits of the possible range  are reached or exceeded as saturation. Interestingly, none of the above-mentioned DP control methods has the drawback that they do not consider input saturation. If saturation is not considered when optimizing control parameters, even a “good” control law  from the objective function's point of view in training may significantly degrade control performance in the natural environment. It has been shown in previous studies that the presence of constraints causes closed-loop systems to be unstable \cite{Stein2003}. Fig.~\ref{fig:Example of Actuator Satuation} shows shows an example of control performance degradation due to actuator saturation. In this example, we assume a situation where the values related to the actuators exceed the constraints. In this case, the actual actuator value given to the control law is replaced by the upper or lower bounds of the constraint, so a significant control output is always obtained. The output value is then used in the next time step to calculate the actuator input value, which also causes the actuator constraint to be exceeded in the next time step. Because of the repetition of these situations, it becomes difficult to resolve the control overshoot, causing the ship to keep missing the target point.

In the study of DPS, Perez et al. proposed a DP control considering disturbances and input saturation using the antiwindup technique, and \cite{Perez2009}. Veksler et al. proposed to deal with actuator saturation as a model predictive control \cite{Veksler2016}. Along with the drive constraints of the actuator, there are also constraints on the time-varying rate of the drive. In this paper, we simply refer to the time-varying rate of the actuator as the rate. In the aerospace field, rate saturation has been studied for many years because it can lead to serious accidents \cite{Hess1997}. In the other hand,ship field , however, only a few studies on the existence of rate are known. For example, Lyngstadaas et al. designed a control model that considers rate saturation \cite{Lyngstadaas2018}. Jin et al. designed an auxiliary dynamical system to deal with rate saturation \cite{Jin2018}.  Suyama et al. addressed the constraint by converting the constrained tracking problem to an unconstrained strict-feedback form \cite{Suyama2023}. In addition, many of such studies are  often mathematically complex.   
\par
Therefore, in this study, we considered the constraints of both actuators and rates only in the parameter optimization stage of the control law. In other words, we use the control law without considering the constraints so as to reduce the cost of the theory construction stage. We   believe that this approach will simplify the design of control laws for  real problems with constraints. In addition, a  path  is needed for the ship to follow   in the automatic navigation. In this paper, we refer to such a path as a reference path . In many previous studies on optimizing control laws, the reference paths required for such tracking control have been predefined. For example, they have been obtained in advance as a time function of a particular form \cite{Do2004}\cite{Dai2016} or generated by a filter based on the current position and endpoint information \cite{Alfheim2018}. In this study, we optimized the parameters in the filter simultaneously with the aforementioned control gains to ensure that the reference generation is applies to  training and testing
  \par
The main contributions of this paper are as follows:
 \begin{enumerate}
    \item We considered optimizing the control gains by setting upper and lower bounds on the actuator control input and rate values.
   \item We considered optimizing the reference filter parameters and control gains to automate the generation of reference paths.
    \item With the optimization described in 1 and the reference partitioning algorithm described in 2, we aimed to obtain references for various paths in a single optimization calculation.
    \end{enumerate}

Section~\ref{sec:notation} presents the different notations used in this study. Section~\ref{sec:System description} describes the coordinate system and dynamics used in this study and explains the backstepping control law and the driving range of the actuators based on the coordinate system and dynamics. Section~\ref{sec:optimizeation} describes the optimization targets for gain and reference generation and the design of evaluation functions to evaluate them. Section~\ref{sec:scenario} describes setting up scenarios for the optimization learning and post-training tests. Section~\ref{sec:simulation result} presents the simulation results. Section~\ref{sec:discussion} discusses our analysis. Finally, Section ~\ref{sec:conclusion} presents the conclusion. 

\section{Notation}\label{sec:notation}
In the following, 1D and 2D represent the one and two dimensional systems, respectively. $\mathbb{R}$ represents the set of real numbers, and $\mathbb{R}^n$ denotes the $n$-dimensional Euclidean space.  The time derivative of a variable $x$ is denoted by $\dot{x}$. The letter ``t'' denotes time.  

\section{System description}\label{sec:System description}
\subsection{Coordinate system}
\begin{figure}[tb]
        \centering 
        \includegraphics[width=1\hsize]{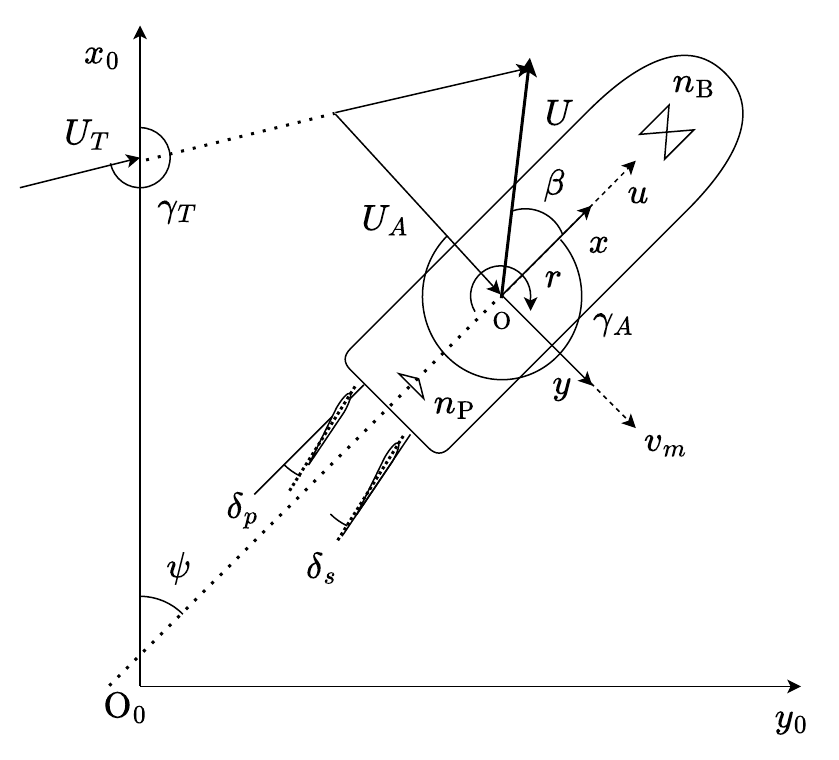}
        \caption{The coordinate systems in this study }
        \label{fig:The Coordinate Systems}
    \end{figure}

 The target ship in this study  is a geometrical similarity model of a ship with a VecTwin rudder, which is a special rudder.
 Fig.~\ref{fig:The Coordinate Systems} shows the coordinate systems of the ship.
 The position of the center of gravity and the heading angle of the ship are expressed in $\mathrm{O}-x_0y_0$ coordinates and denoted by the vector $p(t):= (x_0(t),y_0(t),\psi_0(t))^\top$. 
 \par
 The velocity and Yaw angular velocity of the hull motion are denoted by the vector $v(t):= (v_1(t),v_{\mathrm{m}2}(t),v_3(t))^\top$.
 
\subsection{State equation}
The equation of state is formulated as follows\cite{Sorensen1996}:
 \begin{align}
    \dot{p} &= J(p)v,   \label{eq:coordinate_change} \\
    M\dot{v} + Dv &= \tau + \tau_{\mathrm{wind}} ,\label{eq:force}
 \end{align}
  where $J(p)\in \mathbb{R}^{3\times3}$ is the transformation from hull fixed coordinates to earth fixed coordinates, and 
$M \in \mathbb{R}^{3\times3}$ is the mass and added mass matrix, $D \in \mathbb{R}^{3\times3}$ is the constant decay coefficient matrix, and 
$\tau\;:= (\tau_1,\tau_2,\tau_3)^\top $is the sum of the forces and moments generated by the actuators, and
$\tau_{\mathrm{wind}}\; := 
(\tau_{{\mathrm{wind}},1},\tau_{{\mathrm{wind}},2},\tau_{{\mathrm{wind}},3})^\top$ is the sum of forces and moments due to wind. \\
The matrix in Equations (\ref{eq:coordinate_change})(\ref{eq:force}) are as follows:  

    \begin{eqnarray}
         J(p) :=
            \begin{pmatrix}
            \cos{\psi} & -\sin{\psi} & 0 \\
            \sin{\psi} & \cos{\psi} & 0 \\
            0 & 0 & 1
            \end{pmatrix} 
        \end{eqnarray}
    \begin{eqnarray}
         M :=
            \begin{pmatrix}
            m_{11} & 0 & 0 \\
            0 & m_{22} & m_{23} \\
            0 & m_{32} & m_{33}
            \end{pmatrix} 
        \end{eqnarray}
    \begin{eqnarray}
         D :=
            \begin{pmatrix}
            d_{11} & 0 & 0 \\
            0 & d_{22} & d_{23} \\
            0 & d_{32} & d_{33}
            \end{pmatrix} 
            .
    \end{eqnarray}
  By assuming $|M|\neq 0\;$, Equation (\ref{eq:force}) can be expressed as follows.
    \begin{eqnarray}
         \dot{v} = -M^{-1}Dv + M^{-1}
         (\tau + \tau_{\mathrm{wind}}) \label{eq:-MDv+Mtau}
    \end{eqnarray}
    
    The parameter $\tau$ is calculated as a linear composite of the forces as follows:
\begin{eqnarray}
    \tau = TV\Tilde{u}.
\end{eqnarray}

Also, we use the relation between the actuator input $\Tilde{u}$ and total control force in the VecTwin rudder, as represented by Rachman et al.~\cite{Rachman2023}.  
Equation (\ref{eq:-MDv+Mtau}) is transformed as follows:

\begin{align}
          \dot{v} &= -M^{-1} ( D v + TV\Tilde{u} +   \tau_{\mathrm{wind}} ), 
          \label{eq:-M(Dv+TVu)}\\
           \Tilde{u} &:= (\Tilde{\delta}_\mathrm{P} ,\Tilde{\delta}_\mathrm{S} , \Tilde{n}_\mathrm{B} ) ^\top , \nonumber
\end{align}
where, $\Tilde{\delta}_\mathrm{P} := {\delta}_\mathrm{P} -{\delta}_{\mathrm{P},\mathrm{h}}$, and $\Tilde{\delta}_\mathrm{S} := {\delta}_\mathrm{S} -{\delta}_{\mathrm{S},\mathrm{h}} 
({\delta}_{\mathrm{P},\mathrm{h}},{\delta}_{\mathrm{S},\mathrm{h}} $ are deviation from hover angle)
,  
$\Tilde{n}_\mathrm{B} := n_\mathrm{B} |n_\mathrm{B}|$   is the product of the propeller speed and its absolute value. 


The sum of the forces and moments induced by the wind, $\tau_{\mathrm{wind}}$, is expressed as follows:

\begin{align}
    \tau_{\mathrm{wind}} = 
\begin{bmatrix}
    0.5 \rho_\mathrm{a} |U_\mathrm{A}|^2
    A_{\mathrm{T}} C_{\mathrm{w}x} \\
     0.5 \rho_\mathrm{a} |U_\mathrm{A}|^2
    A_{\mathrm{L}} C_{\mathrm{w}y} \\
     0.5 \rho_\mathrm{a} |U_\mathrm{A}|^2
    A_{\mathrm{L}} L_\mathrm{pp} C_{\mathrm{w}\psi}
\end{bmatrix},   
\end{align}
where, 
       $\rho$ denotes the air density, 
    $L_\mathrm{pp}$ denotes the total length of the ship, 
    $A_{\mathrm{T}},A_{\mathrm{L}}$ denote the cross-sectional area of the ship, 
    $U_\mathrm{A}$ denotes the relative wind speed, 
    $\gamma_\mathrm{A}$ denotes the relative wind direction.
\par

The dimensionless wind coefficient is expressed as follows\cite{ISHERWOOD1973}:
\begin{equation}
\begin{gathered}
    \left\{
\begin{aligned}
   C_{\mathrm{w}x} &= XX_0 + XX_1 \mathrm{cos}( 2\pi - \gamma_\mathrm{A} ) \\
                & \quad+ XX_3 \mathrm{cos}3( 2\pi - \gamma_\mathrm{A} ) + XX_5 \mathrm{cos}5( 2\pi - \gamma_\mathrm{A} ) \\
   C_{\mathrm{w}y} &= YY_1 \mathrm{sin}( 2\pi - \gamma_\mathrm{A} ) \\
    & \quad+ YY_3 \mathrm{sin}3( 2\pi - \gamma_\mathrm{A} ) + YY_5 \mathrm{sin}5( 2\pi - \gamma_\mathrm{A} ) \\
   C_{\mathrm{w}\psi} &= NN_1 \mathrm{sin}( 2\pi - \gamma_\mathrm{A} ) \\
    & \quad + NN_2 \mathrm{sin}2( 2\pi - \gamma_\mathrm{A} ) + NN_3 \mathrm{sin}3( 2\pi - \gamma_\mathrm{A} ).
\end{aligned}
\right.
\end{gathered}
\end{equation}

\subsection{Control law}
A position-following controller was designed based on the backstepping control.
The backstepping method is a typical nonlinear control method to construct a control law for the entire system \cite{krstic1995nonlinear}. 
For a second-order cascade system in the tracking control of a moving object such as a ship, it is possible to design a control law that makes the error system exponentially stable at the origin, and many applications can be found in the field of ship control \cite{Grovlen1996}, \cite{Witkowska2007}.

The control input $\Tilde{u_c}:= (\Tilde{\delta_\mathrm{P}}_c ,\Tilde{\delta_\mathrm{S}}_c , \Tilde{n_\mathrm{B}}_c ) ^\top $ can be chosen as follows. 
For detailed derivation and stability, please refer to the Appendix.

\begin{equation}
    \begin{aligned}
        \Tilde{u}_c &= -B^{-1} J^{-1}(\psi)
     (C_2e_2 +e_1 - C_1^2e_1 + C_1e_2 \\
   &\quad +(\frac{\mathrm{d}J(p)}{\mathrm{d}\psi} {v_3} + J(p)A )v 
   - \ddot{p_\mathrm{d}} - J(p) M^{-1}  \tau_{\mathrm{wind}}). 
   \label{eq:control_wind}
    \end{aligned}
\end{equation}

In the simulations conducted in this study, the actuator state $u:=({\delta}_\mathrm{P} ,{\delta}_\mathrm{S} , {n}_\mathrm{B} ) ^\top $ and its time rate of change are given upper and lower limits. The stern propeller speed is assumed to be constant at $n_P = 60.0~\mathrm{s}^{-1} $. \par

  Considering the mechanical constraints of the model ship and the region where the linear relationship between the VecTwin force and the VecTwin rudder angle is reasonable, the constraints of the actuator state $u$ are defined as follows:
\begin{equation}
    \begin{gathered}
        u \in \mathcal{U} := \mathcal{D}_\mathrm{P} \times \mathcal{D}_\mathrm{S} \times \mathcal{N}_\mathrm{B}\\
        \text{where} \quad \left\{
        \begin{aligned}
            \mathcal{D}_\mathrm{P} &:= [ -\frac{105}{180}\pi,  -\frac{60}{180}\pi ] \\
            \mathcal{D}_\mathrm{S} &:= [ \frac{60}{180}\pi,  \frac{105}{180}\pi] \\
            \mathcal{N}_\mathrm{B} &:= [ -60,  60] .
        \end{aligned}
        \right.
     \end{gathered}
     \label{eq:actuator_limit}
 \end{equation}

 $\Omega := (\Omega_\delta,  \Omega_\delta,  \Omega_B )^\top$ 
 denotes the constant time rate of change of the rudder and bow thruster of the model ship environment.
 Similarly, the time rate of change of the actuator state of the model ship is defined as follows:
     \begin{equation}
         \Omega_\delta = 0.349~\mathrm{s}^{-1}, \; \Omega_\mathrm{B} = 100.0~\mathrm{s}^{-1}.
     \end{equation}

\section{Optimization target and evaluation function}\label{sec:optimizeation}

\subsection{Optimization target}
For parameter optimization, the problem to be solved is formulated as a minimization problem. 
The objective function is defined by the control performance and saturation of the control when a certain candidate solution of the parameters to be optimized is applied to the gain and filter parameters and followed for several scenarios. The scenarios are described in detail in Section~\ref{sec:scenario}.

The objective function $J:\mathbb{R}^{12} \rightarrow \mathbb{R}$
is designed by taking the linear combination of  the indices $J_\mathrm{e}, J_{\widehat{u_\mathrm{c}}}, J_{\widehat{\Delta u}} $ as follows:
\begin{align}
&&J:=  J_\mathrm{e} + r_1 J_{\widehat{u_\mathrm{c}}} + r_2 J_{\widehat{\Delta u}}, 
\end{align}
 where $r_1, r_2$ are weight coefficients.
The following evaluation terms were established for the control input value $u_\mathrm{c}$,
 where $K \in \mathbb{R}$ is the total number of scenarios and $T \in \mathbb{R}$ is the total number of time steps per scenario.
 
 \begin{itemize}
    \item Tracking errors$ : J_\mathrm{e} $
         \begin{eqnarray}
         J_\mathrm{e} := 
          \sum_{k=1}^K
          \sum_{t=1}^T \{ e^{\top}_{k,i} R_1 e_{k,i} \} ,
     \end{eqnarray}
        where,
        \begin{align}
        e &:= 
        (p - p_\mathrm{r}) + w_\mathrm{e} (p - p_\mathrm{d} )
         ,\label{eq:error}
         \end{align}
        
        $R_1 \in \mathbb{R} ^ {3\times3} $
        is the weight matrix of the tracking error and $w_\mathrm{e}$ is the weight coefficient with respect to the reference filter.
     The value of $R_1$ is determined from the investigation of parallel points satisfying Equations (\ref{eq:force}) and (\ref{eq:actuator_limit}) , and the value of $w_\mathrm{e}$ is determined from trial and error as follows.
      \begin{align}
            R_1 &:= \text{diag}\left(  
            \frac{1}{1.0^2}, 
            \frac{1}{{1.0}^2}, 
            \frac{1}{(0.2{\pi})^2}\right), \\
            w_\mathrm{e} &:= 10.
        \end{align}
     
  \item Exceeding saturation of actuator state$ : J_{\widehat{u_\mathrm{c}}} $
     \begin{equation}
         J_{\widehat{{u_\mathrm{c}}}} :=
          \sum_{k=1}^K 
           \sum_{i=1}^T\{\widehat{{u_\mathrm{c}}}^{\top}_{k, i} R_2 {\widehat{{u_\mathrm{c}}}}_{k, i} \}.
     \end{equation}
     
     ${\widehat{u_\mathrm{c}} } \in \mathbb{R}^3 $
     is the amount by which the control value $u_\mathrm{c}$ exceeds its upper and lower bounds.\\
     Each component of ${\widehat{u_c}}_j (j = 1,2,3 )$ is defined by the following equation.
      \begin{eqnarray}
         {\widehat{u_c}}_j := {u_c}_j - \mathrm{clip}({u_c}_j, {u_c}_j^{\mathrm{min}},{u_c}_j^{\mathrm{max}} ), \label{eq:input_rate}
    \end{eqnarray}

    where, the function $\mathrm{clip}$ is a function that places the signal $s\in R $ in the interval $[s_\mathrm{min}, s_{\mathrm{max}}]$.
    \begin{align}
    \mathrm{clip}(s,s_{\mathrm{min}}, s_{\mathrm{max}}) = 
        \begin{cases} 
            s_{\mathrm{max}}\; &\mathrm{for} \;s_{\mathrm{max}} \geq s \\
            s \; &\mathrm{for} \;s_{\mathrm{min}} \geq s > s_{\mathrm{max}} \\
            s_{\mathrm{min}} \;&\mathrm{for}\; s < s_{\mathrm{min}} .
        \end{cases}
    \end{align}
         $R_2 \in \mathbb{R}^{3 \times 3}$
     is the penalty weight matrix for $\widehat{u_\mathrm{c}}$ and is defined as follows.
    \begin{equation}
        R_2 := \frac{r_1}{3} \mathrm{diag}\left( \frac{1}{\Delta {\delta_\mathrm{P}}^2},  \frac{1}{\Delta {\delta_\mathrm{S}}^2},  \frac{1}{\Delta {n_\mathrm{B}}^2} \right) ,
    \end{equation}
    \begin{align}
        \begin{cases}
           \Delta {\delta_\mathrm{P}} &:= \mathrm{max}\mathcal{D}_\mathrm{P} - \mathrm{min}\mathcal{D}_\mathrm{P} \\
           \Delta {\delta_\mathrm{S}} &:= \mathrm{max}\mathcal{D}_\mathrm{S} - \mathrm{min}\mathcal{D}_\mathrm{S} \\
           \Delta {n_\mathrm{B}} &:= \mathrm{max}\mathcal{N}_\mathrm{B} - \mathrm{min}\mathcal{N}_\mathrm{B} .
        \end{cases}
    \end{align}

   \item Exceeding rate saturation of actuator state:$ J_{\widehat{\Delta u}} $
    \begin{equation}
        J_{\widehat{\Delta u}} :=
        \sum_{k=1}^K   \sum_{i=1}^T \{\widehat{\Delta u}^{\top}_{k, i} R_3{\widehat{\Delta u}}_{k, i} \},
    \end{equation}
    where $\widehat{\Delta u} \in \mathbb{R} ^ 3$
    is the amount by which the time variation of the transition from state $u$ to the control value $u_\mathrm{c}$ during $\Delta t$ exceeds its upper and lower bounds.
    \\
    Each component $\Delta u_j (j = 1,2,3 )$ of $\Delta u$ is defined by the following equation.
    \begin{align}
        \widehat{\Delta u_j} &:= \Delta u_j - \mathrm{clip}(\Delta u_j - \Omega_j,\Omega_j ),\\
        \Delta u &:= \frac{ u_c - u }{\Delta t}.
    \end{align}
    $R_3 \in \mathbb{R}^{3\times3}$
    is the penalty weight matrix for $\widehat{\Delta u}$ and is defined as follows :
    \begin{equation}
        R_3 := \frac{r_2}{3} \text{diag}\left(\frac{1}{{\Omega _\delta}^2},  \frac{1}{{\Omega _\delta}^2},  \frac{1}{{\Omega _\mathrm{B}}^2} \right).
    \end{equation}
    \end{itemize}

\subsubsection{Control law gain}
We optimize the design matrix contained in the control law, as shown in Equation(\ref{eq:control_wind}). 
For positive symmetric matrices $C_1, C_2$, the Cholesky decomposition is expressed as follows:
\begin{equation}
    C_1 = A_1 A_1{^\top},~ C_2 = A_2 A_2^{\top},
\end{equation}
where,  $A_1$ and $A_2$ are the lower triangular matrices expressed as follows:
\begin{align}
     A_1 = 
    \begin{bmatrix}
                    a_{11} & 0  & 0 \\
                    a_{12} & a_{13}  & 0 \\
                    a_{14} & a_{15}  & a_{16} \\
    \end{bmatrix}, 
    A_2 = 
    \begin{bmatrix}
                    a_{21} & 0  & 0 \\
                    a_{22} & a_{23}  & 0 \\
                    a_{24} & a_{25}  & a_{26} \\
    \end{bmatrix}  .
\end{align}

From the above, we optimize 12 independent parameters of $C_1, C_2$: 
\begin{equation}
    X_\mathrm{C} = (a_{11},  \dots , a_{16}, a_{21},  \dots , a_{26}) \in \mathbb{R}^{12} .\label{eq:X_C_parameter}
\end{equation}
\par
For the search range of each parameter, considering the condition that $C_1, C_2$ are positive definite and the accuracy, we set
The search range for each parameter is defined as follows, considering the condition that $C_1, C_2$ are positive definite.

\begin{equation}
    \begin{gathered}
         \left\{
        \begin{aligned}
0.001 &\leq a_{11}, a_{13}, a_{16}, a_{21}, a_{23}, a_{26} \leq 10 \\
-10 &\leq a_{12}, a_{14}, a_{15}, a_{22}, a_{24}, a_{25} \leq 10.
        \end{aligned}
        \right.
     \end{gathered}
     \label{eq:search_range}
 \end{equation}

\subsubsection{Parameters in reference filters}
 The reference filter used in this study was proposed in the book by Fossen \cite{Fossen2011}. In \cite{Alfheim2018}, the filter was used to generate the reference paths for the proposed DP control system follows.
 \par
 It consists of a first-order low-pass filter.
 By applying a filter to the reference input value $p_\mathrm{r} = [x_\mathrm{r},y_\mathrm{r},\psi_\mathrm{r}]^\top$, the target position and attitude $p_\mathrm{d} = [x_\mathrm{d},y_\mathrm{d},\psi_\mathrm{d}]^\top$ is obtained by applying a filter to the target position and attitude $p_\mathrm{d}$ .

 \begin{align}
    &\dddot{p}_\mathrm{d} + (2\Delta + I_3)\Omega \ddot{p}_\mathrm{d} 
    + (2\Delta + I_3) {\Omega}^2 \dot{p}_\mathrm{d} +{\Omega}^3 p_\mathrm{d} 
    = {\Omega}^3 p_\mathrm{r} \nonumber \\
    \Leftrightarrow  &\dddot{p}_\mathrm{d} 
    + G \Omega \ddot{p}_\mathrm{d} 
    + G {\Omega}^2 \dot{p}_\mathrm{d} +{\Omega}^3 p_\mathrm{d} 
    = {\Omega}^3 p_\mathrm{r} \label{filter},
\end{align}   
where, 
    $\Delta = \mathrm{diag}[\zeta_x,\zeta_y,\zeta_{\psi}]^\top,  
    \zeta$ denotes the relative damping ratio, 
    $\Omega = \mathrm{diag}[\omega_{nx},\omega_{ny},\omega_{n\psi}]^\top,
    \omega$ denotes the natural frequency, 
    and $G = 2\Delta + I_3$.
\par
From the above,
we optimize six independent parameters of  $\Delta,\Omega $ :
\begin{equation}
    X_{\mathrm{ref}} = (\zeta_x,\zeta_y,\zeta_{\psi},
    \omega_{nx},\omega_{ny},\omega_{n\psi}
    ) \in \mathbb{R}^6 .\label{eq:X_ref_param}
\end{equation}

The search ranges are defined as follows:
\begin{equation}
    \begin{gathered}
         \left\{
        \begin{aligned}
   0.01 &\leq \zeta_x,\zeta_y,\zeta_{\psi} 
 \leq  0.1\\
  0.8 &\leq \omega_{nx},\omega_{ny},\omega_{n\psi} \leq 2.0~.
        \end{aligned}\right.
    \end{gathered}
\end{equation}

From Equations (\ref{eq:X_C_parameter}) and (\ref{eq:X_ref_param}), for a total of 18 independent parameters, CMA-ES (Covariance Matrix Adaptation Evolution Strategy). It is an  excellent black-box optimization algorithm that introduces variance-covariance into the evolutionary strategy\cite{Hansen2016}.
Compared to other algorithms\cite{Tomera2014}\cite{Witkowska2012}, CMA-ES with the introduction of covariance is the most efficient for solbing optimization problems with relatively complex search geometries. It has also been used for parameter tuning in PID control\cite{Henclova2019}.
\section{Scenario planning}\label{sec:scenario}
In designing scenarios, it is important to align the state quantities of values such as position and velocity during training and testing. This is because the performance of the controller cannot be accurately measured. If they are not aligned, the results of training will not be reflected in the test, and as a result, the performance may not necessarily improve.
\subsection{Training scenarios}
The time in one scenario is defined $t = 120~\mathrm{s}$, and the simulation time range is  $\mathrm{d}t = 0.1~\mathrm{s}$.
During the training scenarios, we assumed that the wind  blows from all eight directions in the simulations. The wind speed is assumed to be constant at $U = 1.0~\mathrm{m/s}$.
Table~\ref{tab:object} shows the reference input value $p_\mathrm{r}$ of each episode.
\begin{table}
  \begin{center}
  \caption{Target position and attitude for each scenario}
 \begin{tabular}{|c|c|c|c|} \hline
   Episode & $x~[\mathrm{m}]$ & $y~[\mathrm{m}]$ & $\psi~[\mathrm{degree}]$ \\ \hline\hline
   1 & 0 & 0 & 0 \\
   2 & 4.0 & 0 & 0 \\ 
   3 & $-$4.0 & 0 & 0 \\
   4 & 0 & 4.0 & 0 \\
   5 & 0 & $-$4.0 & 0 \\
   6 & 0 & 0 & 30 \\
   7 & 0 & 0 & $-$30 \\ 
   8 & 4.0 & 4.0 & 0 \\ 
   9 & 4.0 & 0 & 30 \\ 
   10 & 0 & 4.0 & 30 \\ 
   11 & 4.0 & 4.0 & 30 \\ \hline
 \end{tabular}
 \label{tab:object}
 \end{center}
 \end{table}
 
In this study, we optimized  the control law without considering input and rate saturation is also performed. Comparing the results with and without saturation shows that our proposed approach performs well enough in practice
\par 
We varied the weight coefficients in the cases with and without saturation.  Table~\ref{tab:weight} shows the weights of the objective function in each case.
\begin{table}
  \begin{center}
  \caption{Weights of the objective function used in cases}
 \begin{tabular}{|c|c|c|c|} \hline
   Case  & $r_1$ & $r_2$ \\ \hline\hline
   1  & 10 & 10 \\
   2  & 0 & 0 \\ \hline
 \end{tabular}
 \label{tab:weight}
 \end{center}
 \end{table}
 
\subsection{Test scenarios}
We conducted a test scenario based on the 4-corner DP test presented in \cite{Skjetne2017} to evaluate the performance of the optimized control law. In this scenario, the four corners start from the origin, which is one of the corners, followed in sequence and return to the origin. Therefore, in a single trial, it is possible to evaluate the movements in the sway, surge, and yaw directions and the combined movements in the sway-yaw and surge-yaw directions.
\par
The detailed behavior of the test scenario is as follows :
\begin{enumerate}
\item Change $+5~\mathrm{m}$ position in the surge direction.
\item Change $+5~\mathrm{m}$ position in the sway direction.
\item Change direction $+\frac{45}{180}\pi$ in the yaw direction.
\item Change $-5~\mathrm{m}$ position in the sway direction while maintaining the course.
\item Change $-\frac{45}{180}\pi$ in the yaw direction and $-5~\mathrm{m}$ in the surge direction.
\end{enumerate}
\begin{figure}[tb]
        \centering 
        \includegraphics[width=1\hsize]{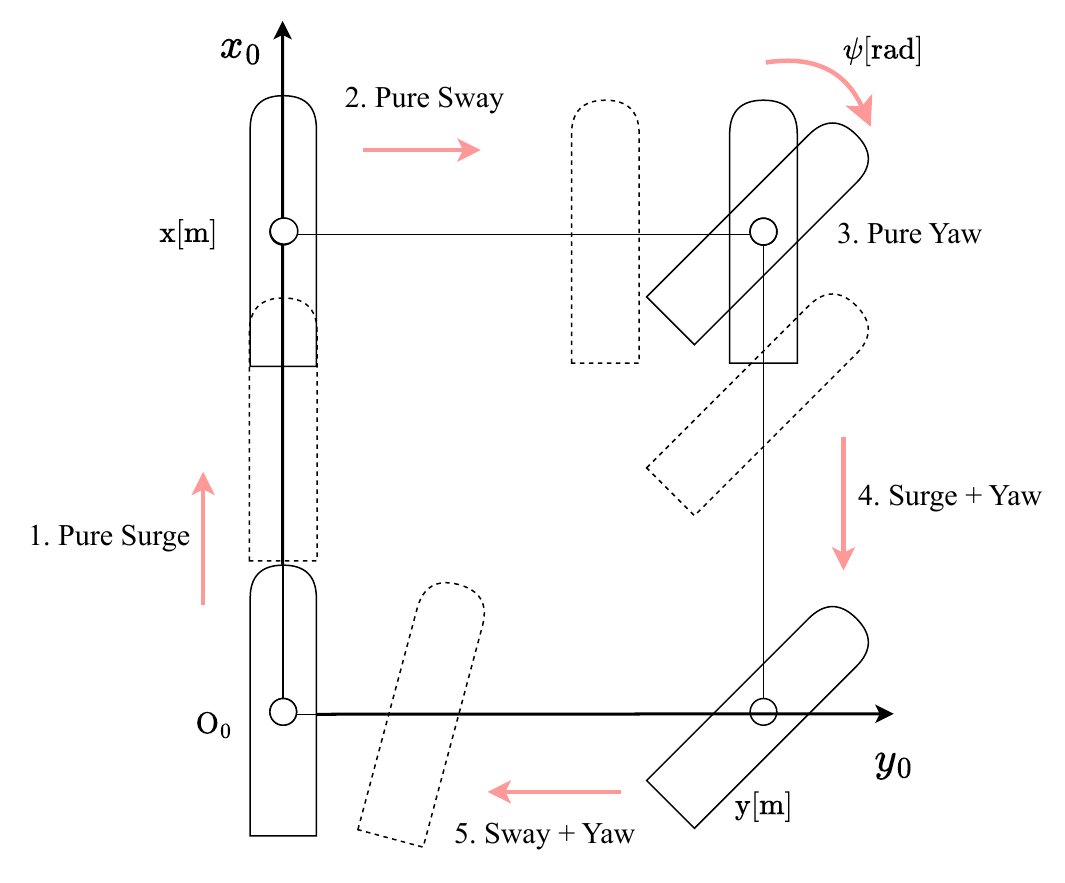}
        \caption{Example of 4 corner testing}
        \label{fig:four_corner}
    \end{figure}
In the test scenario, the wind speed is assumed to be constant at $U = 0.5~\mathrm{m/s}$, and the wind direction is assumed to be constant at $\gamma = \frac{30}{180}\pi$  .

\subsection{Target position and attitude segmentation algorithm}
Even if an appropriate target is set, e.g., too large, an extreme target position setting may make the control system unstable. In this study, we designed an algorithm in which references at given at regular intervals to achieve a highly scalable reference setting. The target positions and attitudes  in the optimization calculations are set according to the interval of the partitioning.  This is expected to alleviate the instability of control for large targets and improve the versatility of the control system.
\section{Simulation result}\label{sec:simulation result}
\begin{figure}[tb]
        \centering 
        \includegraphics[width=1.2\hsize]{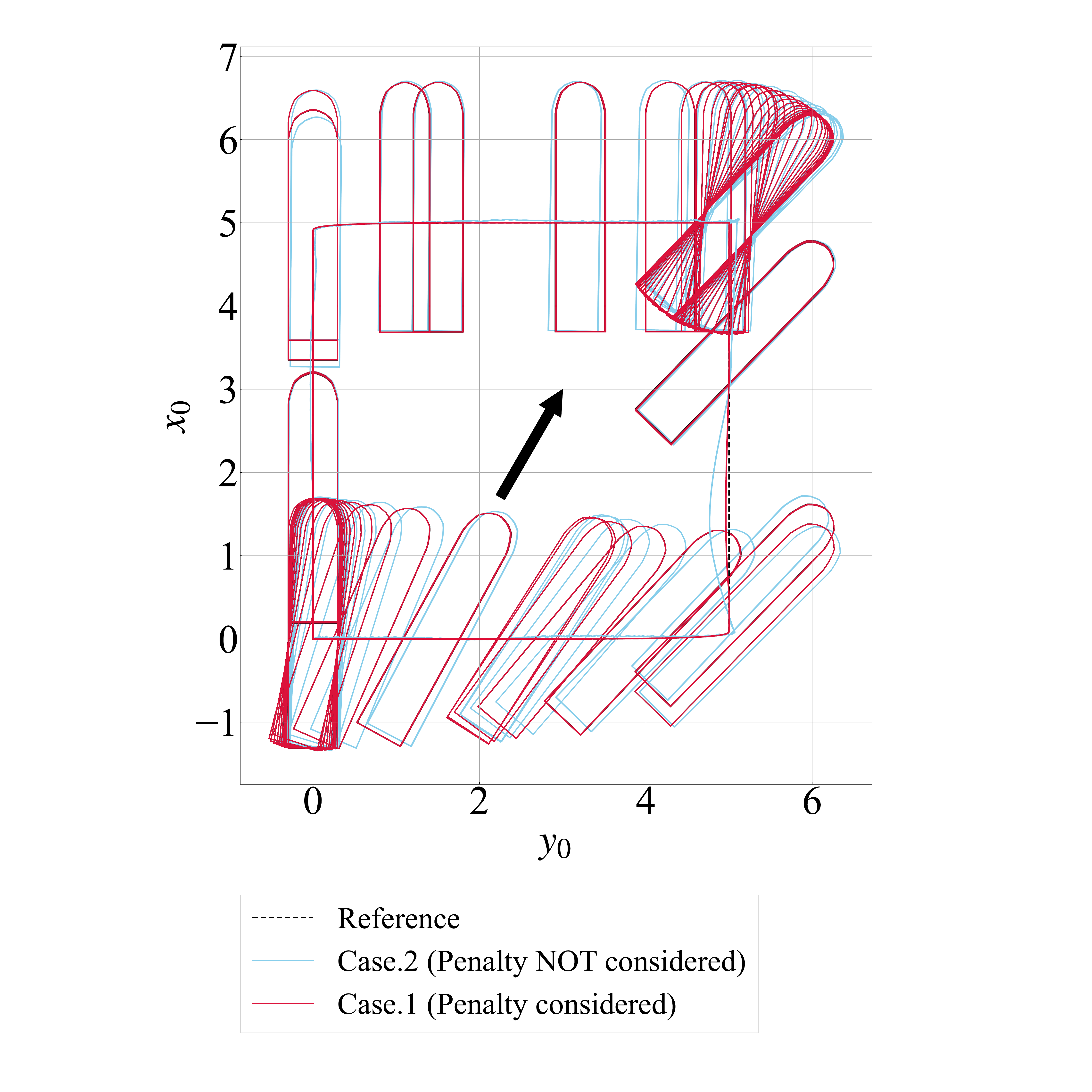}
        \caption{The results of the trajectory for the test scenario.
        }
        \label{fig:Result of Trajectory}
    \end{figure}
\begin{figure*}
    \centering
\includegraphics[width=1\hsize]{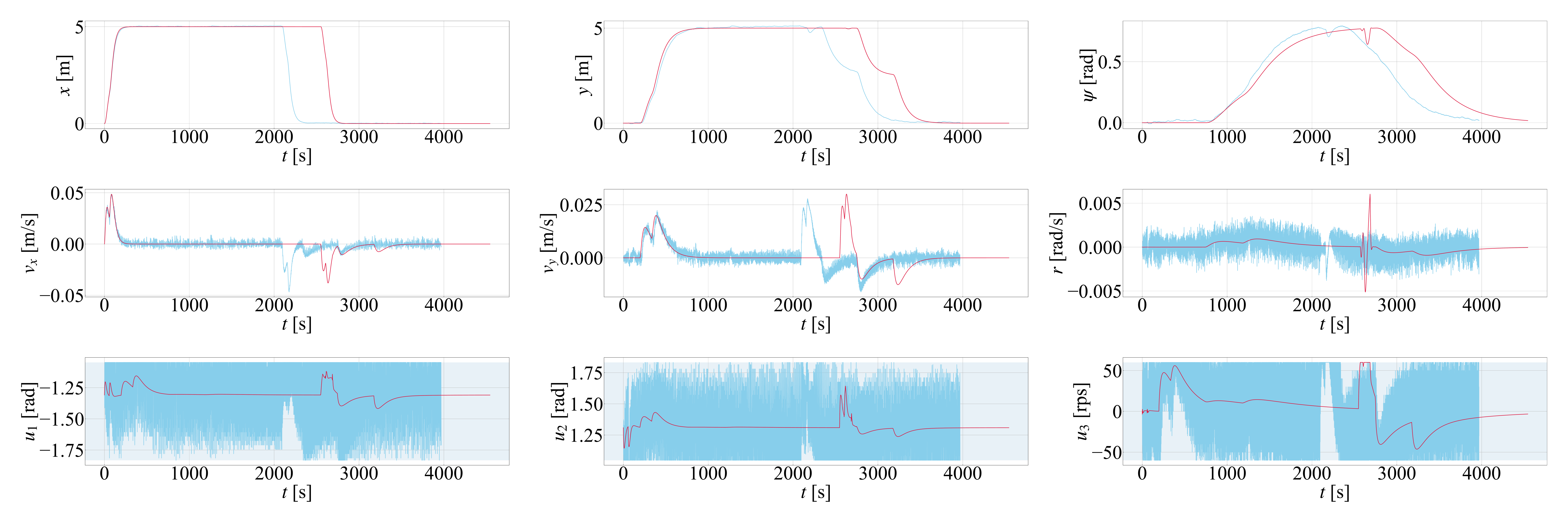}
    \caption{
     Time variation of each parameter related to the control law in the test scenario. 
     The parameters $u_1,u_2,$ and $u_3$ represent $\delta_{\mathrm{P}},\delta_{\mathrm{S}},$ and $n_\mathrm{B}$, respectively. 
     The light blue areas indicate the driving range. The colors of the lines in the graph are the same as in the legend of Fig.~\ref{fig:Result of Trajectory}}
    
    \label{fig:Result of Parameters}
\end{figure*}

Figs.~\ref{fig:Result of Trajectory} and ~\ref{fig:Result of Parameters} show the test scenario results using the parameter values optimized by the training scenario. The parameters obtained in Case. 1 (red line), where both penalties for input and rate saturation were considered, can follow the reference path. On the other hand, Case. 2 (blue line), where both penalties were not considered, shows that the behavior deviates from the reference path.

The reference (black line) generated by the filter, which is optimized simultaneously with the control parameters, is considered to generate the appropriate distance target position and posture for the behavior by the controller.
\section{Discussion and limitation}\label{sec:discussion}
We treated the actuator drive input and rate constraints as penalties and examined saturation to control gain optimization. As a result, we obtained control gains that do not degrade the tracking performance  under the disturbance of wind conditions while minimizing value exceedances. Also, we obtained parameter values that generate a suitable reference path for such control.
\par
 Now, we discuss the challenges and future developments in this study. The biggest challenge in our designed control law is its weakness against coupled motion combining multiple movements. The simulation results show that the actuators are slightly exceeded in a coupled two-directional motion (see Fig.~\ref{fig:Result of Parameters}). We believe this may be due to a lack of learning for complex coupled motions and insufficient power output of the thrusters of the assumed vessels. Therefore, we propose further research into the causes and solutions. In this study, we determined the values of some hyperparameters by trial and error. For example, the weight coefficients in the evaluation function equations for the target position and tracking error in the test scenarios were determined manually. The relationship between these values and the values to be optimized is still unclear, and this connection needs to be established. In addition, the disturbance was limited to the effect of steady wind. To simulate a ship in an environment closer to the actual operating environment, it is necessary to consider control under other complex disturbance conditions, such as irregular winds and disturbances caused by water waves. One of the future tasks for the authors is to determine the subfishing force of irregular winds using the fast calculation method of wind disturbance presented by Maki et al. \cite{Maki2022}.
\section{Conclusion}\label{sec:conclusion}\label{sec:conclusion}
%
We have studied the treatment of constraints on state, control input, and their rates in automatic control algorithms for ships. We designed a control law for the DPS control based on the backstepping method and optimized its control parameters using CMA-ES. The parameters of the reference generation filter were also optimized to enable the conversion to various paths. The simulation results showed that control gains can be obtained with less constraint exceedance and appropriate reference paths can be generated. 

\section*{Conflict of interest}

    The authors declare that they have no conflict of interest.

\bibliographystyle{spphys}       
\bibliography{reference}   

\section*{Appendix(Detailed derivation and the stability analysis of the control law)}
    \subsection{Detailed derivation of the control law}
    Define the error $e_1$ between the current position and attitude$p:= [x_0,y_0,\psi_0]^\top$ and the target position and attitude $p_\mathrm{d} = [x_\mathrm{d},y_\mathrm{d},\psi_\mathrm{d}]^\top$as follows~\cite{Krikelis1984}:
    \begin{equation}
        e_1 := p - p_\mathrm{d} 
    \end{equation}
     For the time derivative of $e_1$, we have:
     \begin{equation}
         \dot{e_1} = J(p)v - \dot{p_\mathrm{d}}
     \end{equation}
     The virtual control $\alpha_1$ is considered to consist of a second backstepping variable and a stabilizing function, given by:
     \begin{equation}
         \alpha_1 = -C_1e_1 + \dot{p_\mathrm{d}}
     \end{equation}
     where, $C_1 \in \mathbb{R}^{3\times3}$is a strictly positive definite design matrix.\par
     Denote $e_2$ as follows:
     \begin{equation}
         e_2 := \dot{p} - \alpha_1
     \end{equation}
     Here, $\dot{e_1}$ can be expressed as:
     \begin{equation}
         \dot{e_1} = -{C_1} {e_1} + e_2 \label{eq:e_1}
     \end{equation}
     The time derivative of $e_2$ is expressed as follows:
     \begin{align}
         \dot{e_2} &= \ddot{p} - \dot{\alpha_1} \nonumber\\
    &= \frac{d}{dt}(J(p)v)+C_1\dot{e_1} - \ddot{p_\mathrm{d}} \nonumber\\
    &=C_1 \dot{e_1} + \frac{dJ(p)}{d\psi} v_3 v + J(p)\dot{v} - \ddot{p_\mathrm{d}} \nonumber \\
    &= C_1(-C_1e_1 + e_2 ) + \frac{dJ(p)}{d\psi} v_3 v + J(p)(Av + B\Tilde{u} + \tau_\mathrm{wind} ) - \ddot{p_\mathrm{d}}\nonumber \\ 
    &= \left(\frac{dJ(p)}{d\psi} v_3 - J(p)A\right)v 
     +J(p)B\Tilde{u} \\ &~~~~~~~+ C_1\left(-C_1e_1 + e_2 \right)- J(p)M^{-1}\tau_\mathrm{wind} -\ddot{p_\mathrm{d}}
     \end{align}
     where, $A:= M^{-1}D,~ B:= M^{-1}TV \in \mathbb{R}^{3\times3}$\par
     The control input $\Tilde{u^c}:= (\Tilde{\delta^c}_\mathrm{P} ,\Tilde{\delta^c}_\mathrm{S} , \Tilde{n^c}_\mathrm{B} ) ^\top $ can be chosen as in Equation (\ref{eq:control_wind}).
     Where, assume that $C_2\in\mathbb{R}^{3\times3}$ can be expressed as follows:
     \begin{equation}
         \dot{e}_2 = -{C_2}{e_2} - {e_1}\label{eq:e_2}
     \end{equation}
     
     \subsection{Stability analysis of the control law}
 From Equations (\ref{eq:e_1}) and (\ref{eq:e_2}), we obtain as follows:
\begin{equation}
       \dot{z} = -Cz + Sz  
\end{equation}
where:
\begin{align}
    z &:= ({e_1}^\top,{e_2}^\top) ^\top \in \mathbb{R}^6 \\
    C &:= \mathrm{diag}(C_1, C_2) \in \mathbb{R}^{6\times6}\\
    S &:= \begin{pmatrix}
        O & I_3 \\
        -I_3 & O \\
    \end{pmatrix}
    \in \mathbb{R}^{6\times6}
\end{align}

At this point, we can define a candidate Lyapunov function as follows:
\begin{equation}
    V := \frac{1}{2} z^\top z > 0 , {\forall}z \neq 0.\label{eq:lyapunov}
\end{equation}

Differentiating Equation (\ref{eq:lyapunov})  with respect to z over time, we obtain the following relation
\begin{align}
    \dot{V} &= z^{\top} \dot{z} ~= z^{\top} (-Cz + Sz) \nonumber \\
    &= -z^{\top} Cz < 0 , {\forall}z \neq 0 \label{eq:dot_lyapunov}
\end{align}
From Euation (\ref{eq:dot_lyapunov}), the control law is stable.

\end{document}